\documentclass{amsart}
\usepackage{amssymb}
\usepackage{amsfonts}
\usepackage{amsmath}
\usepackage{lastpage}
\usepackage[legalpaper,bookmarks=true,colorlinks=true,linkcolor=blue,citecolor=blue]{hyperref}
\usepackage{graphicx}%
\usepackage{fancyhdr}
\usepackage{color}
\usepackage[mathlines]{lineno}
\usepackage{lscape}
\usepackage{epsfig}
\usepackage{natbib}
\usepackage[]{listings}
\usepackage{geometry}
\usepackage{tgbonum}
\fontfamily{qcr}\selectfont

\newtheorem{theorem}{Theorem}
\theoremstyle{plain}

\newtheorem{corollary}{Corollary}

\newtheorem{lemma}{Lemma}

\numberwithin{equation}{section}

\begin{document}
\Large

\title[Flexible versions of the Stone-Weierstrass Theorem]{Flexible versions of the Stone-Weierstrass Theorem in General and Applications to Probability Theory}
\author{Gane Samb Lo}

\begin{abstract} When applying the classical Stone-Weierstrass common version in Probability Theory for example, and in other fields as well, problems may arise if all points of the compact set are not separated. A solution may consist in going back to the proof and finding alternative versions. In this note, we did it and come back with two flexible versions which are easily used for the needs of classical Probability Theory.\\

\noindent  Gane Samb Lo.\\
LERSTAD, Gaston Berger University, Saint-Louis, S\'en\'egal (main affiliation).\newline
LSTA, Pierre and Marie Curie University, Paris VI, France.\newline
AUST - African University of Sciences and Technology, Abuja, Nigeria\\
gane-samb.lo@edu.ugb.sn, gslo@aust.edu.ng, ganesamblo@ganesamblo.net\\
Permanent address : 1178 Evanston Dr NW T3P 0J9,Calgary, Alberta, Canada.\\

\end{abstract}
\maketitle

\section{Introduction}

\noindent In Probability Theory, the theorem of Stone-Weierstrass is one of the possible tools to prove that the characteristic function serves as an effective characterization of the probability law of a random vector (say, in $\mathbb{R}^k$, $k\geq 1$) and that the weak convergence of Probability laws is equivalent to the convergence of characteristics functions.\\

\noindent However, the application of such a theorem from the common version of the theorem is not straightforward. An appeal for an extended version (see  \cite{simmons}, page 165) is made. An example is available in  \cite{billingsley}, but such a version is stated on locally Hausdorff Compact sets and uses functions vanishing at infinity. At a such an earlier stage of Probability Theory, it is regrettable that a quite highly sophisticated approach is required. Beyond this example, we think that the simple knowledge of the common version only does not allow the researcher in Probability and Statistics to draw all the benefits of that extraordinary tool. More generally, it seems that authors work from the consequences rather than the true principle.\\

\noindent Here we want to show that simple versions can give more flexibility in applying that Theorem by highlighting the elements of the proofs and by sticking more to the lines of the proof. In doing so, we just allow flexibility in the application of the mentioned theorem.\\

\noindent The rest of the paper is organized as follows. In the Section \ref{sec2}, we just state the most common version of the Stone-Weierstrass theorem stated in the introductory books of Topology. In the same section, we explain the reasons that led us to the current extended version. In Section \ref{sec3}, we provide remarks from the proof of the theorem in classical textbooks and just propose versions preserving the classical conclusions. In Section \ref{sec4}, we show how to prove the Billingsley problem with the version of this paper.

\section{Common versions of the Stone Weierstrass Theorem and its application} \label{sec2}
    
The theorem is stated as follows.

\begin{theorem} (Complex version of Stone-Weierstrass's Theorem) \label{sec_EF_theo_02} Let $K$ be a non-singleton compactum (A Hausdorff space on which the Heine-Borel holds) and $\mathcal{A}$ be a non-empty sub-algebra of $C(K,\mathbb{C})$, the continuous functions defined on $K$ with values in ring $\mathbb{C}$ of complex numbers, such that the following assertion :\\

\noindent (1) A separates the points of $K$, that is, for any distinct elements $x$ and $y$ of $K$, there exists $f \in A$ such that $f(x) \neq f(y)$,\\

\noindent and one of the two assertions \\

\noindent (2) For any $x \in K$, there exists $f \in \mathcal{A}$ such that $f(x) \neq 0$.\\

\noindent (3) $A$ contains all the constant functions,\\

\noindent and the additional assertion :\\

\noindent (4) For all $f \in \mathcal{A}$, its conjugate function $\bar{f}=\mathcal{R}(f) - i \mathcal{I}m(f) \in A$,\\

\noindent  hold. Then $\mathcal{A}$ is dense in $C(K,\mathbb{C})$
$$
\overline{A}=C(K,\mathbb{C}).
$$
\end{theorem}

\noindent The case where $K$ is a singleton is dismissed as an obvious thing.\\

\noindent This theorem is applied in Probability Theory and in many other fields like Neural networks XXXX. Let us explain the following case study. Recall that two probability measures $\mathbb{P}_1$ and $\mathbb{P}_2$ on a metric space $(E,d)$ endowed with the Borel $\sigma$-algebra $\mathbb{B}(E)$ \textit{are equal} if and only if : for any $f$ in the class $C_b(E)$ of bounded and continuous real-valued functions defined on $E$, we have

$$
\int_E f d\mathbb{P}_1=\int_E f d\mathbb{P}_2. \ \ (EL)
$$ 

\bigskip \noindent (See \cite{ips-mfpt-ang}, Chapter 3, Part III for a proof). If $E=\mathbb{R}^k$, we define for any probability measure $\mathbb{P}$ on $(\mathbb{R}^k, \mathbb{B}(\mathbb{R}^k)$, its characteristic function by

$$
\Psi_{\mathbb{P}}(u)=\int_{\mathbb{R}^k} \exp(i \langle u, x\rangle) \mathbb{P}(x), \ u \in \mathbb{R}^k.
$$

\noindent \noindent As in classical books of Probability Theory, we want to show that the equality of the characteristic functions 

$$
\Psi_{\mathbb{P}_1}=\Psi_{\mathbb{P}_1}, \ \ (EC)
$$ 

\noindent on $\mathbb{R}^k$ is equivalent to Formula (EL) and hence to the equality in law. In fact, that (EL) implies (EC) is obvious. The main work is to prove that (EC) implies (EL). The main ideas in the proof are :\\

\noindent (a) Fixing a real number $a>0$ and considering $K_a=[-a,a]^k$.\\

\noindent (b) Considering the family $\mathcal{H}$ of functions which are finite linear combination (with real coefficients) of the form

$$
\exp(i\pi \langle n, u\rangle /a), \ u\in \mathbb{R}^k.
$$

\bigskip \noindent where $n=(n_1,...,n_k)^T \in \mathbb{Z}^k$.\\

\noindent (c) Showing that the class of $\mathcal{H}_a$ of restrictions of elements of $\mathcal{H}$ on $K_a$ is dense in $C(K_a,\mathbb{R})$ by the Stone-Weierstrass theorem, which implies that for any $\varepsilon>0$, for any $f \in C(\mathbb{R}^k)$, there exists $h \in \mathcal{H}$ such that

$$
\|f-h\|_{K_a}=\sup_{x \in K_a} |f(x)-h(x)|\leq \varepsilon,
$$

\bigskip \noindent The two following properties are important in the sequel :\\

\noindent (d) By periodicity, the uniform norm of $h$ on $K_a$ is equal to the uniform norm of $h$ on the whole space.\\

\noindent (e) The integral of $h$ with respect to  $\mathbb{P}_i$, on $\mathbb{R}^k$, is a finite linear combination of the values of the characteristic function of $\mathbb{P}_i$, $i \in \{1,2\}$. 

\noindent The conclusion is achieved through a smart combination of limit results as $a \uparrow +\infty$. But the key tool is in Point (c).\\

\noindent \textbf{Unfortunately}, the subclass $\mathcal{H}_a$ possesses all the desired conditions to apply Stone-Weierstrass's theorem on the compact $K_a$ except the separation of the points of $K_a$. Indeed, each $h \in \mathcal{H}$ assigns the same values of all $2^k$ edge points of $k_a$ of the form $(\pm a, \pm a, \cdots, \pm a)^T$.\\

\noindent The version we are going to provide will solve this drawback.\\

\noindent Before we go further, it is worth-mentioning that in some special and particular cases, it might be possible to show that $\overline{A}=C(K,\mathbb{C})$ by explicitly constructing for any $f\in C(K,\mathbb{C})$ a sequence $(f_n)_{n\geq 0} \subset \mathcal{A}$ such that $f_n$ converges uniformly to $f$ as $n$ goes to infinity. But it is not reasonable to expect this is non simple cases.\\

\section{Analysis of the Proof of Stone-Weierstrass's theorem} \label{sec3}

\noindent The proof as in classical textbooks is based on the two following lemmas.

\begin{lemma} \label{sec_EF_lem_01} Let $K$ be a compactum (A Hausdorff space on which the Heine-Borel holds) and $\mathcal{A}$ be a non-empty lattice subclass of $C(K)$, that is $\mathcal{A}$ is closed under finite minimum and maximum of its elements :

$$
\forall (f,g)\in A, \ min(f,g) \in \mathcal{A} \ and \ max(f,g) \in \mathcal{A}.
$$

\bigskip \noindent Let $f\in C(K)$. Then $f \in \overline{\mathcal{A}}$ if and only if for each two distinct elements $x$ and $y$ of $K$, the restriction $f_{|\{x,y\}}$ is limit of a restrictions of  a sequence of elements of $\mathcal{A}$ on $\{x,y\}$ :\\

$$
\forall (x,y) \in K^2, \ \exists (f_{n}^{x,y})_{n\geq 0} \subset \mathcal{A}, \ f_{n}^{x,y}(t) \rightarrow f(t) \ as \ n\rightarrow +\infty, \ for \ t \in \{x,y\}. \label{C2a2}
$$ 
\end{lemma}

\begin{lemma} \label{sec_EF_lem_02} Let $\mathcal{A}$ be a  closed non-empty sub-algebra in $C(K)$. Then $\mathcal{A}$ is a lattice space.\\

\noindent As a consequence, a closed non-empty sub-algebra in $C(K)$ is lattice.
\end{lemma}

\noindent By combining these two lemmas, proving $\overline{\mathcal{A}}=C(K,\mathbb{C})$ is the same as proving that $\overline{\mathcal{A}_1}=C(K,\mathbb{C})$ with  
$\mathcal{A}_1=\overline{\mathcal{A}}$, which is a lattice sub-algebra of $C(K,\mathbb{C})$ whenever $\mathcal{A}$ is a sub-algebra of $C(K,\mathbb{C})$. This leads to the following general form

\begin{theorem} (A first general version of Stone-Weierstrass's Theorem) \label{sec_EF_theo_02} Let $K$ be a compact space (A Hausdorff space on which the Heine-Borel holds) and $\mathcal{A}$ be a non-empty sub-algebra of $C(K,\mathbb{C})$, the continuous functions defined on $K$ with values in the ring $\mathbb{C}$ of complex numbers. We have $\overline{\mathcal{A}}=C(K,\mathbb{C})$ whenever we have for all $f \in C(K,\mathbb{C})$,

\begin{equation}
\forall (x,y) \in K^2, \ \exists (f_{n}^{x,y})_{n\geq 0} \subset \mathcal{A}, \ f_{n}^{x,y}(t) \rightarrow f(t) as n\rightarrow +\infty, for \ t \in \{x,y\}. \label{C2a2a}
\end{equation} 
\end{theorem}

\bigskip \noindent By stating this, we simply bring to the surface the most inner technical tool in the proof of the classical Theorem.\\

\noindent Actually the usual version of the Stone-Weierstrass version seeks to get this by requiring Assumptions (1) and (2) or (4) in Theorem \ref{sec_EF_theo_02} in the real version. Assumption (4) is required to extend the real version to the complex version. Let us remain in the real case. The proof in \cite{choquet1966} and in almost many other books goes too far. Let us explain why.\\

\noindent Let $x$ and $y$ be two distinct elements of $K$. Let $f \in C(K)$ with $f(x)=\alpha$ and $f(y)=\beta$. In the classical proof, we combine Assumptions (1) and (2) to find a function $g \in A$ such that $g(x)=\alpha$ and $g(y)=\beta$. So, Formula \ref{C2a2a} holds but with a \textbf{constant sequence} $f_{n}^{x,y}=f^{x,y}$. If $x=y$, and since $K$ is not a singleton, the method is re-conducted for $x$ and $z$ with $x \neq z$ and Formula \ref{C2a2a} also holds.\\

\noindent \textbf{As a conclusion}, making happen Formula \ref{C2a2a} with a constant sequence is a high price to pay for getting the approximation. The following version of the Stone-Weierstrass Theorem will be still valid. It is based on the fact that we do not need that $\mathcal{A}$ separates all the points. 

%\begin{theorem} (A second general version of Stone-Weierstrass's Theorem) \label{sec_EF_theo_02} Let $K$ be a compact space (A Hausdorff space on which the Heine-Borel holds) and $\mathcal{A}$ be a non-empty sub-algebra of $C(K,\mathbb{C})$, the continuous functions defined on $K$ with values in the ring $\mathbb{C}$ of complex numbers. We have $\overline{\mathcal{A}}=C(K,\mathbb{C})$ whenever we have

%$$
%\forall (x,y) \in K^2, \ \exists (f_{n}^{x,y})_{n\geq 0} \subset \mathcal{A}, \ f_{n}^{x,y}(t) \rightarrow f(t) \ as \ n\rightarrow +\infty, \ for \ t \in \{x,y\}. \label{C2a2a}
%$$ 
%\end{theorem}

\begin{corollary} (A second version of Stone-Weierstrass's Theorem) \label{sec_EF_theo_04} Let $K$ be a non-singleton compactum (A Hausdorff space on which the Heine-Borel holds) and $\mathcal{A}$ be a non-empty sub-algebra of $C(K,\mathbb{C})$, the continuous functions defined on $K$ with values in ring $\mathbb{C}$ of complex numbers. Suppose that there exists $K_0 \subset K$ such that we have

$$
\forall (x,y) \in K_0^2, \ \exists (f_{n}^{x,y})_{n\geq 0} \subset \mathcal{A}, \ f_{n}^{x,y}(t) \rightarrow f(t) \ as \ n\rightarrow +\infty, \ for \ t \in \{x,y\}. \label{C2a2b}
$$ 

\noindent (1) $\mathcal{A}$ separates the points of $K\setminus K_0$ and separates any point of $K_0$ from any point of $K$, \\

\noindent and one of the two assertions \\

\noindent (2) For any $x \in K\setminus K_0$, there exists $f \in \mathcal{A}$ such that $f(x) \neq 0$.\\

\noindent (3) $A$ contains all the constant functions,\\

\noindent and the additional assertion :\\

\noindent (4) For all $f \in \mathcal{A}$, its conjugate function $\bar{f}=\mathcal{R}(f) - i \mathcal{I}m(f) \in A$,\\

\noindent  hold. Then $\mathcal{A}$ is dense in $C(K,\mathbb{C})$ :
$$
\overline{A}=C(K,\mathbb{C}).
$$
\end{corollary}

\begin{corollary} (A third version of Stone-Weierstrass's Theorem) \label{sec_EF_theo_05} Let $K$ be a non-singleton compactum (A Hausdorff space on which the Heine-Borel holds) and $\mathcal{A}$ be a non-empty sub-algebra of $C(K,\mathbb{C})$, the continuous functions defined on $K$ with values in ring $\mathbb{C}$ of complex numbers. Suppose that there exists $K_0 \subset K$ such that $K\setminus K_0$ has at least two elements and the following assertions :

\noindent (0) We have : for all $f \in C(K,\mathbb{C})$

$$
\exists (f_{n})_{n\geq 0} \subset \mathcal{A}, \ \forall x \in K_0, \  f_{n}(x) \rightarrow f(x) \ as \ n\rightarrow +\infty, \ for \ t \in \{x,y\}. \label{C2a2b}
$$ 

\noindent or \\

$$
\exists g \in \mathcal{A}, \ \forall x \in K_0,  \ \ g(x)=f(x),
$$ 

\bigskip \noindent (1) $\mathcal{A}$ separates the points of $K\setminus K_0$ and separates any point of $K_0$ from any point of $K$, \\

\noindent and one of the two assertions \\

\noindent (2) For any $x \in K\setminus K_0$, there exists $f \in \mathcal{A}$ such that $f(x) \neq 0$.\\

\noindent (3) $\mathcal{A}$ contains all the constant functions,\\

\noindent and the additional assertion :\\

\noindent (4) For all $f \in \mathcal{A}$, its conjugate function $\bar{f}=\mathcal{R}(f) - i \mathcal{I}m(f)$ belongs to $\mathcal{A}$.\\

\noindent  Then $\mathcal{A}$ is dense in $C(K,\mathbb{C})$ :
$$
\overline{A}=C(K,\mathbb{C}).
$$
\end{corollary}

\section{Return to the Application of \textit{SW} Theorem in Probability Theory} \label{sec4}

Let us go back to Section \ref{sec2}, Point (b). We are going to apply Corollary \ref{sec_EF_cor_05}. Assume we have the same notation (regarding $]a,b[$, $r>0$ and $K_r$). Let us denote $C_{b,0,r}(\mathbb{R}^d)$ the restrictions of elements of $C_{b,0}(\mathbb{R}^d)$ on $K_r$. We should prove that $C_{b,0,r} \subset \overline{\mathcal{H}_r}$. For this, let us take $f \in C_{b,0,r}(\mathbb{R}^d)$. The general form of an element of $\mathcal{H}$ is

$$
h(x)=\sum_{1 \leq j \leq p} a_p \exp(i\pi\langle n_p, x\rangle/r), \ x \in \mathbb{R}^k,
$$ 

\noindent where $a_p \in \mathbb{C}$,  $n_p \in \mathbb{Z}^k$. Let us take

$$
K_0=\partial K_r=\{x \in K_r : \forall \ j \in \{1,\cdots,k\}, \ x_i=-r \ or \ x_i=r\}.
$$

\noindent We have $f=0$ on $K_0$. If $x$ and $y$ are two points in $K_r$ such they are not among the edging points in the border $\partial K_r$ both, then there exists $j_0 \in \{1,\cdots,d\}$ such that $0<|x_{j_0}-y_{j_0}|<2r$ that is $|(x_{j_0}-y_{j_0})/r|<2$ and the function $h_{r}(x)=\exp(i\pi x_{j_0}/r)$ separates $x$ and $y$. By adding the other assumptions, conclude that $f$ is in the closure of $\mathcal{H}_r$.\\

\end{document}